\documentclass[12pt]{amsart}

\usepackage[cp1251]{inputenc}
\usepackage[english]{babel}

\usepackage{vmargin}
\setpapersize{A4}
\setmarginsrb{2.5cm}{0cm}{2.5cm}{1cm}{2cm}{1cm}{0cm}{1.5cm}

\usepackage{amsmath, amsfonts, amssymb, amsthm}
\usepackage{rotating, multirow}
\usepackage{comment}

\newcommand\geqs{\geq}

\newcommand\q[1]{\mathfrak{q}(#1)}
\newcommand\hc{\chi}

\newtheorem{prop}{Proposition}
\newtheorem{cons}{Corollary}

\theoremstyle{definition}
\newtheorem*{prf}{Proof}

\title[An analogue of the Perelomov-Popov formula]
{An analogue of the Perelomov-Popov formula
\\[4pt] 
for the Lie superalgebra $\q{N}$}

\author{T. A. Grigoryev}
\address{
T. A. Grigoryev: Skolkovo Institute of Science and Technology, Skolkovo Innovation Center, 
Moscow 143026, Russia and 
National Research University Higher School of Economics, Department of Mathematics, 
Moscow 119048, Russia}

\author{M. L. Nazarov}
\address{M. L. Nazarov:
Department of Mathematics, University of York, York YO10 5DD, United Kingdom}

\begin{document}

\maketitle

\thispagestyle{empty}

\begin{abstract}
We study the center of the universal enveloping algebra of the strange Lie 
superalgebra $\q{N}$. We obtain an analogue of the well known Perelomov-Popov formula \cite{popov} for central elements of this algebra -- an
expression of the central characters through the highest weight parameters.
\end{abstract}

\section{Introduction}
\subsection{Lie superalgebra $\q{N}$}

The strange Lie superalgebra $\q{N}$ can be realized as a subalgebra in the
general linear Lie superalgebra $\mathfrak{gl}(N|N)$ over the complex field,
see for instance \cite{cheng_wang}. If the indices $i$ and $j$ 
range over $-N, \ldots, -1, 1, \ldots N$ then 
the $4N^2$ elements $E_{ij}$ 
span the algebra $\mathfrak{gl}(N|N)$ as a vector space. 
The Lie superbracket on $\mathfrak{gl}(N|N)$
is defined~by
\begin{equation*}
[E_{ij}, E_{kl}] = \delta_{jk} E_{il} - (-1)^{(\bar{\imath} + \bar{\jmath})(\bar{k} + \bar{l})}\delta_{il} E_{kj}
\end{equation*}
where 
$$
\overline{k} = \begin{cases} 
\,0
\quad\text{if}&k > 0\,,
\\ 
\,1
\quad\text{if}&k < 0\,.
\end{cases}
$$
Then $\q{N} \subset \mathfrak{gl}(N|N)$ is the subalgebra
of fixed points of the involution
$$\eta: E_{ij} \mapsto E_{-i-j}\,.$$
Hence as a vector space $\q{N}$ is spanned by the $2N^2$ elements
$$F_{ij} = E_{ij} + E_{-i-j}$$
with $i > 0$. The Lie superbracket on $\q{N}$ is then described by
\begin{equation}
\label{commut}
[F_{ij}, F_{kl}] = 
\delta_{kj} F_{il} - (-1)^{(\overline{\imath}+ \overline{\jmath})(\overline{k} + \overline{l})} \delta_{il} F_{kj} + \delta_{k-j} F_{-il} - (-1)^{(\overline{\imath} + \overline{\jmath})(\overline{k} + \overline{l})} \delta_{-il} F_{k-j}\,.
\end{equation}

\subsection{Casimir elements in $U(\q{N})$}

Unless otherwise stated, we will be assuming that 
the indices in the expressions below range over
$-N, \ldots, -1,$ $ 1, \ldots, N$.
For any $n=1,2,\ldots$ consider 
the elements of the universal enveloping algebra $U(\q{N})$
first proposed in \cite{sergeev1}:
\begin{equation}
\label{casimirs}
C^{(n)}_{ij} = 
\sum_{k_1,\ldots,k_{n-1}}F_{ik_1} (-1)^{\overline{k}_1} F_{k_1k_2} (-1)^{\overline{k}_2} \cdots F_{k_{n-2}k_{n-1}} (-1)^{\overline{k}_{n-1}} F_{k_{n-1}j}\,.
\end{equation}
Taking into account the recurrence relations
\begin{equation*}
C^{(n+1)}_{ij} = 
\sum_k
F_{ik}\, (-1)^{\overline{k}}\, C^{(n)}_{kj}
\end{equation*}
one can find the supercommutator
\begin{equation*}
[F_{ij}, C^{(n)}_{kl}]= 
\delta_{kj} C^{(n)}_{il}-(-1)^{(\overline{\imath} + \overline{\jmath})(\overline{k} + \overline{l})} \delta_{il} C^{(n)}_{kj} + \delta_{k-j} 
C^{(n)}_{-il} - (-1)^{(\overline{\imath} + \overline{\jmath})
(\overline{k} + \overline{l})} \delta_{-il} C^{(n)}_{k-j}
\end{equation*}
which is similar to the superbracket (\ref{commut}) between the
generators of $\q{N}\,$. It is then easy to see that the elements 
$$
c_n = \sum_{i} C^{(n)}_{ii}
$$ are central in $U(\q{N})\,$. Moreover, by using the relation
\begin{equation*}
C^{(n)}_{-i-j} = (-1)^{n-1}\, C^{(n)}_{ij}
\end{equation*}
following from (\ref{casimirs}) we see that 
$c_n=0$ if $n$ is even. This is why we will be only interested in
the $c_n$ with $n$ odd. These are the {\it Casimir elements} for the 
Lie superalgebra $\q{N}$ as introduced in \cite{sergeev2}.
It was shown 
in \cite{nazarov_sergeev} that these elements generate 
the center of $U(\q{N})\,$.

\subsection{Harish-Chandra homomorphism}

Let $v$ be a singular vector of any irreducible
finite-dimensional representation of the Lie superalgebra $\q{N}$ relative to the natural triangular decomposition 
$\q{N} = \mathfrak{n}_{-} \oplus \mathfrak{h} \oplus \mathfrak{n}_{+}\,$ where
\begin{align*}
\begin{split}
\mathfrak{n}_{-} & = {\rm span}\left\{F_{ij} \,\big|\, |i| > |j|\right\},\\
\mathfrak{h} & = {\rm span}\left\{F_{ij} \,\big|\, |i| = |j|\right\},\\
\mathfrak{n}_{+} & = {\rm span}\left\{F_{ij} \,\big|\, |i| < |j|\right\}.
\end{split}
\end{align*}
Then the following equalities hold:
\begin{align}
\label{highestVec1}
F_{ij} \cdot v &= 0 
\quad\text{for}\quad 
|i| < |j|\,,
\\
\label{highestVec2}
F_{ii} \cdot v &= \lambda_i\,v
\quad\text{for}\quad
i>0\,.
\end{align}
Here the $\lambda_i$ are the eigenvalues of the elements $F_{ii} = F_{-i-i}$ of the even part of the Cartan subalgebra $\mathfrak{h}_0 = {\rm span}\left\{F_{ii}\, \big|\, i>0\right\}$. 
They depend on the particular representation of $\q{N}\,$.
Let $\lambda \in \mathfrak{h}_0^*$ be the highest weight of the
representation so that $\lambda(F_{ii})=\lambda_i$ for $i>0\,$.

The generators of $\q{N}$ in the summands of (\ref{casimirs}) can always be rearranged in such a way that in each monomial left-to-right first go the lowering operators, then the operators from Cartan subalgebra and last the raising operators. 
Here the lowering operators are elements of $\mathfrak{n}_-$ and the raising operators are elements of $\mathfrak{n}_+\,$. 
The operators from Cartan subalgebra can also be rearranged so 
that elements from its even part $\mathfrak{h}_0$ go after those from its odd part $\mathfrak{h}_1 = {\rm span}\left\{F_{-ii}\, \big|\, i>0\right\}\,$. 
It suffices to use the supercommutation relations (\ref{commut}) to achieve this. 
The part of the resulting sum which belongs to $U(\mathfrak{h})$
is well-defined. For the sum (\ref{casimirs}) with $i=j$
this is its image
under the {\it Harish-Chandra homomorphism} $\chi\,$, see~\cite{cheng_wang}.

The subspace formed by the vectors of an irreducible representation of $\q{N}$
satisfying (\ref{highestVec1}),(\ref{highestVec2}) 
is called the singular subspace. The peculiarity of $\q{N}$ 
shows in the fact that the singular subspace of an irreducible representation 
is not one-dimensional, but is irreducible over the Cartan subalgebra 
$\mathfrak{h} = \mathfrak{h}_0 \oplus \mathfrak{h}_1$ for which we have
$U(\mathfrak{h}) = S(\mathfrak{h}_0) \otimes \wedge(\mathfrak{h}_1)\,$,
see \cite{P}.
 
Due to (\ref{casimirs}) the Casimir elements are even and therefore 
commute with the whole algebra $U(\q{N})$ in the usual non 
$\mathbb{Z}_2$-graded sense. 
This implies that they act as scalar operators in the irreducible representations. This allows us to consider their eigenvalues when acting on some fixed singular vector instead and the computations in this case are rather simple.

Let $n$ be odd and 
$v$ be a singular vector of an irreducible representation of 
$\q{N}\,$. Then
\begin{equation}
\label{action}
c_n\cdot v = \hc(c_n)\cdot v = 
\hc(c_n)\,\big|_{\,F_{ii} = \lambda_i} \, v
\quad\text{for}\quad
i>0\,.
\end{equation}
For $i>0$ we will automatically substitute
$F_{ii} \mapsto \lambda_i$ after applying
the 
homomorphism $\hc$ as we did in (\ref{action}).
Hence we will be describing the action on the singular vector explicitly.


\section{Computations}

\subsection{Recurrence relations}

Here we will derive a recurrence relation for
the images of the elements $C^{(n)}_{ii}$ with $n=1,3,\ldots$
under the Harish-Chandra homomorphism. 
For brevity we will denote by $G_i$ the element 
$F_{-ii} = F_{i-i}$ of the odd part of the Cartan~subalgebra.
 
\begin{prop}
For $i>0$ we have $\hc(G_i^2) = \lambda_i\,$.
\end{prop}

\begin{prf}
We have
$G_i^2 \cdot v = F_{-ii}^2 \cdot v = \frac{1}{2}\, [F_{-ii}, F_{-ii}] \cdot v= F_{ii} \cdot v = \lambda_i\,v$
for $i>0\,.$
\end{prf}

\begin{prop}
We have
$C^{(n)}_{ij} \cdot v = 0$ whenever $|i| < |j|\,$.
\end{prop}

\begin{prf}
Suppose that $|i| < |j|$. Then $C^{(1)}_{ij} \cdot v = F_{ij} \cdot v = 0\,$. 
Let us use the induction on $n\,$:
$$
C^{(n+1)}_{ij} \cdot v =\sum_k 
F_{ik}\,(-1)^{\overline{k}}\,C^{(n)}_{kj} \cdot v = \sum_{|k| \geqs |j|} (-1)^{\overline{k}}\,[F_{ik}, C^{(n)}_{kj}] \cdot v = \sum_{|k| \geqs |j|} (-1)^{\overline{k}}\,C^{(n)}_{ij} \cdot v = 0\,.
$$
\end{prf}

\begin{prop}
For $i>0$ we have
$$
\hc(C^{(n+1)}_{ii}) = \lambda_i\, \hc(C^{(n)}_{ii}) - G_i\, \hc(C^{(n)}_{-ii}) - \sum_{|k| > i}\, \hc(C^{(n)}_{kk})\,.
$$
\end{prop} 

\begin{prf}
For $i>0$ the vector $C^{(n+1)}_{ii} \cdot v$ equals
\begin{gather*}
\sum_{|k| \geqs i} F_{ik}\, (-1)^{\overline{k}}\, C^{(n)}_{ki} \cdot v =
(F_{ii} C^{(n)}_{ii} - F_{i-i} C^{(n)}_{-ii}) \cdot v + \sum_{|k| > i} (-1)^{\overline{k}}\, [F_{ik}, C^{(n)}_{ki}] \cdot v=
\\
(\lambda_i C^{(n)}_{ii} - G_i C^{(n)}_{-ii}) \cdot v + \sum_{|k| >i} (-1)^{\overline{k}}\, \big(C^{(n)}_{ii} - 
(-1)^{\overline{k}}\, C^{(n)}_{kk}\big) \cdot v= 
\\
\big(\lambda_i C^{(n)}_{ii} - 
G_i C^{(n)}_{-ii} - \sum_{|k|>i} C^{(n)}_{kk}\big) \cdot v\,. 
\end{gather*}
\end{prf}

\begin{cons}
For $i>0$ and $m=1,2,\ldots$ 
we get $\hc(C^{(2m+1)}_{ii}) = 
\lambda_i\,\hc(C^{(2m)}_{ii}) - G_i\,\hc(C^{(2m)}_{-ii})\,$.
\end{cons}

\begin{prop}
For $i>0$ we have
$$
\hc(C^{(n+1)}_{-ii}) =
G_i \hc(C^{(n)}_{ii}) - \lambda_i \hc(C^{(n)}_{-ii}) - \sum_{|k|>i} (-1)^{\overline{k}}\, \hc(C^{(n)}_{k-k})\,.
$$
\end{prop} 

\begin{prf}
For $i>0$ the vector $C^{(n+1)}_{-ii} \cdot v$ equals
\begin{gather*}
\sum_{|k| \geqs i} F_{-ik}\, 
(-1)^{\overline{k}}\,C^{(n)}_{ki} \cdot v = 
(F_{-ii} C^{(n)}_{ii} - F_{ii} C^{(n)}_{-ii}) \cdot v + \sum_{|k| >i} (-1)^{\overline{k}}\, [F_{-ik}, C^{(n)}_{ki}] \cdot v =
\\
(G_i C^{(n)}_{ii} - \lambda_i C^{(n)}_{-ii}) \cdot v + \sum_{|k| > i} (-1)^{\overline{k}}\,\big(C^{(n)}_{-ii} - 
(-1)^{\overline{k}\,(\overline{k} + 1)}\,C^{(n)}_{k-k}\big) \cdot v =
\\
\big(G_i C^{(n)}_{ii} - \lambda_i C^{(n)}_{-ii}\big) \cdot v - \sum_{|k|>i} (-1)^{\overline{k}}\, C^{(n)}_{k-k} \cdot v\,. 
\end{gather*}
\end{prf}

\begin{cons}
For $i>0$ and $m=1,2,\ldots$ we get 
$\hc(C^{(2m)}_{-ii}) = 
G_i \hc(C^{(2m-1)}_{ii}) - \lambda_i\hc(C^{(2m-1)}_{-ii})\,$.
\end{cons}

\begin{prop}
\label{P5}
For $i>0$ and $m=1,2,\ldots$ we have the relation
$$
\hc(C^{(2m+1)}_{ii}) = \hc(C^{(2m+1)}_{-i-i}) = \lambda_i (\lambda_i - 1) \hc(C^{(2m-1)}_{ii}) - 2\lambda_i \sum_{j>i}^{} \hc(C^{(2m-1)}_{jj})\,.
$$
\end{prop} 

\begin{prf}
If $i>0$ then the vector $C^{(2m+1)}_{ii}\cdot v $ equals
\begin{gather*}
\big(\lambda_i C^{(2m)}_{ii} - G_i C^{(2m)}_{-ii}\big) \cdot v =
\\
\lambda_i \big(\lambda_i\, C^{(2m-1)}_{ii} - G_i\, C^{(2m-1)}_{-ii} - 
\sum_{|j| > i}^{}\, C^{(2m-1)}_{jj}\big) \cdot v - G_i \big(G_i C^{(2m-1)}_{ii} - \lambda_i C^{(2m-1)}_{-ii}\big) \cdot v =
\\
(\lambda_i^2 - G_i^2) C^{(2m-1)}_{ii} \cdot v - \lambda_i \sum_{|j|>i} C^{(2m-1)}_{jj} \cdot v = \lambda_i (\lambda_i - 1) C^{(2m-1)}_{ii} \cdot v - 2\lambda_i \sum_{j>i} C^{(2m-1)}_{jj} \cdot v\,. 
\end{gather*}
\end{prf}

\begin{cons}
For $i>0$ and $m=0,1,2,\ldots$ we have
$\displaystyle
\hc(C^{(2m+1)}_{ii}) = \sum_{j = 1}^N\big( A^{m}\big)_{ij}\,\lambda_j$ where
$$A =
\begin{pmatrix}
\lambda_1(\lambda_1 - 1) & -2\lambda_1 & \cdots & -2\lambda_1\\
0 & \lambda_2(\lambda_2 - 1) & \cdots & -2\lambda_2\\
\vdots & & \ddots & \vdots\\
0 & \cdots & 0 & \lambda_N(\lambda_N - 1)
\end{pmatrix}.
$$
\end{cons}

\begin{prf}
This follows from Proposition \ref{P5} by
taking into account that $\hc(C^{(1)}_{ii}) = \hc(F_{ii}) = \lambda_i\,$.
\end{prf}

\subsection{Generating functions}

In order to compute $\hc(c_{2m+1})$ more explicitly, for each $i > 0$
consider the generating function
\begin{equation}
\label{pi}
\mu_i(u) = \sum_{m=0}^\infty \hc(C^{(2m+1)}_{ii})\,u^{-2m-1} = u \sum_{j=1}^N \big((u^2 - A)^{-1}\big)_{ij} \lambda_j
\end{equation}
and write 
$$
A=\Lambda(\Lambda - 1) - 2\Lambda \Delta (1 - \Delta)^{-1}
$$ where $\Delta_{ij} = \delta_{i\,j-1}$ and $\Lambda = \rm{diag}(\lambda_i)$. Denote
$$
\Pi = \big(u^2 - \Lambda(\Lambda + 1)\big) \big(u^2 - \Lambda(\Lambda - 1)\big)^{-1}\,.
$$ 
Then
\begin{equation*}
(u^2 - A)^{-1} = \frac{1}{2}\, (1 - \Pi)\, (1 - \Delta \Pi)^{-1}\, (1 - \Delta) \Lambda^{-1}.
\end{equation*}
The last two factors here cause all summands but the last in the sum (\ref{pi}) cancel, leaving
\begin{gather*}
\mu_i(u) = \frac{u}{2}\, (1 - \Pi_{ii}) \left((1 - \Delta \Pi)^{-1}\right)_{iN} =\\
\frac{u}{2}\, (1 - \Pi_{ii})\, (\Delta \Pi)_{i\,i+1} (\Delta \Pi)_{i+1\,i+2} \ldots (\Delta \Pi)_{N-1\,N} = \frac{u}{2}\, (1 - \Pi_{ii}) \prod_{j > i} \Pi_{jj}\,.
\end{gather*}
This can also be written in more straightforward way:
\begin{equation*}
\mu_i(u) = \frac{u\lambda_i}{u^2 - \lambda_i(\lambda_i - 1)}
\,\prod_{j > i}\, 
\frac{u^2 - \lambda_j(\lambda_j + 1)}{u^2 - \lambda_j(\lambda_j - 1)}\,.
\end{equation*}
In order to find a generating function $\mu(u)$ of the images of the 
central elements under $\chi\,$, we should now sum all
the above obtained expressions for $\mu_i(u)$ 
over the positive values of the index $i$ and recall that 
$C_{-i-i}^{(2m+1)} = C_{ii}^{(2m+1)}$. By making further cancellations we get
\begin{equation*}
\mu(u) = \sum_{m=0}^\infty \hc(c_{2m+1}) u^{-2m-1}=2\,\sum_{i=1}^N\,\mu_i(u)= u \left(1 - \prod_{j=1}^N \Pi_{jj}\right)\,.
\end{equation*}
Finally, we obtain 
\begin{equation*}
\mu(u) = u \left(1 - \prod_{i=1}^N \frac{u^2 - \lambda_i(\lambda_i + 1)}{u^2 - \lambda_i(\lambda_i - 1)}\right).
\end{equation*}

\subsection{Images of the central elements}

Introduce a new variable $z = u^{-2}$ and define 
\begin{equation*}
\widetilde{\mu}(z) = u\, \mu(u) = \sum_{m=0}^\infty \hc(c_{\,2m+1})\,z^{m} 
=\frac{1}{z} \left(1 - \prod_{i=1}^N \frac{1 - z\lambda_i(\lambda_i + 1)}{1 - z\lambda_i(\lambda_i - 1)}\right).
\end{equation*}
From the above definition it immediately follows that 
$\hc(c_{\,2m+1}) = {\rm Res}_0\; \widetilde{\mu}(z) z^{-m-1}\, dz$. 
Combining this with our explicit expression for $\widetilde{\mu}(z)$
we obtain that $\hc(c_{\,2m+1})$ equals
$$
-\sum_{i=1}^N {\rm Res}_{\left(\lambda_i(\lambda_i - 1)\right)^{-1}}\; \widetilde{\mu}(z)\,z^{-m-1}\, dz = \sum_{i=1}^N {\rm Res}_{\left(\lambda_i(\lambda_i - 1)\right)^{-1}}\, \prod_{j=1}^N\, 
\frac{1 - z\lambda_j(\lambda_j + 1)}{1 - z\lambda_j(\lambda_j - 1)}\, z^{-m-2}\, dz
$$
where the regularity of the form at the infinity is taken into account. 
Finally for $m\geqslant0$
$$\hc(c_{\,2m+1}) = 2\sum_{i=1}^N \lambda_i^{m+1}(\lambda_i - 1)^m \prod_{j\not=i} \frac{\lambda_i(\lambda_i - 1) - \lambda_j(\lambda_j + 1)}{\lambda_i(\lambda_i - 1) - \lambda_j(\lambda_j - 1)}\ .
$$
Note that despite the fact that the last expression is formally a 
rational function of $\lambda\,$, its denominator always cancels
which allows us to regard $\hc(c_{\,2m+1})$ as a polynomial of $\lambda\,$.

\section*{Remarks}
The analogue of the Perelomov-Popov formula for the Lie superalgebra $\q{N}$ presented here was obtained by the second named author about 30 years ago but left unpublished. 
Independently but by the same method, this analogue was obtained in \cite{BK} and recently by the first named author. Publishing this note gives us an opportunity to review the history
of this analogue. We thank Jonathan Brundan, Maria Gorelik, Alexandre Kirillov, Grigori Olshanski, Ivan Penkov, Vera Serganova and Alexander Sergeev for helpful conversations.


\begin{thebibliography}{0}
\addcontentsline{toc}{section}{\bibname}

\bibitem{BK}
J. Brundan and A. Kleshchev,
{\it Modular representations of the supergroup Q(n), I\/},
J. Algebra 
{\bf 260}
(2003),
64--98.

\bibitem{cheng_wang} 
S. Cheng and W. Wang,
{\it Dualities and representations of Lie superalgebras},
AMS, Providence, 2013.

\bibitem{sergeev2} 
D. Leites and A. Sergeev,
{\it  Casimir operators for Lie superalgebras},
in: E. Ivanov {\it et al\/} (eds.),
\lq\lq Supersymmetries and Quantum Symmetries\rq\rq,
JINR, Dubna, 2000, pp. 409--411.

\bibitem{nazarov_sergeev}
M. Nazarov and A. Sergeev,
{\it Centralizer construction of the Yangian of the queer Lie superalgebra},
in: J. Bernstein {\it et al\/} (eds.),
\lq\lq Studies in Lie Theory\rq\rq,
Boston, Birkha\"user, 2006,
pp. 417--441. 


\bibitem{P}
{I. B. Penkov},
{\it Characters of typical irreducible finite-dimensional
${\mathfrak q}(n)$-modules},
Funct. Anal. Appl.   
{\bf 20}
(1986),
30--37.

\bibitem{popov} 
A. M. Perelomov and V. S. Popov, 
{\it Casimir operators for semisimple Lie groups},
Math. USSR Izv.
{\bf 2} (1968),
1313--1335.
 
\bibitem{sergeev1} 
A. Sergeev,
{\it The centre of enveloping algebra for Lie superalgebra Q(n,$\mathbb{C}$)}, Lett. Math. Phys.
{\bf 7} (1983),
177--179.


\end{thebibliography}
\end{document}